\pdfoutput=1
\documentclass[a4paper,12pt]{amsart}
\usepackage{amsmath,amssymb,amsthm}
\usepackage[english]{babel}
\usepackage[T1]{fontenc}
\usepackage{newtxtext,newtxmath}
\usepackage[utf8]{inputenc}
\usepackage[babel=true]{microtype}
\usepackage[autostyle=true]{csquotes}
\usepackage{mathrsfs}
\usepackage{thmtools,thm-restate}
\usepackage[pdfusetitle]{hyperref}
\usepackage[citestyle=numeric,bibstyle=numeric,backend=biber,url=false,doi=true,isbn=false,giveninits]{biblatex}

\AtEveryBibitem{\clearfield{month}}
\AtEveryBibitem{\clearfield{day}}
\AtEveryBibitem{\clearfield{language}}

\bibliography{literature.bib}
\usepackage{tikz-cd}
\usepackage{tikz}
\usetikzlibrary{calc,3d}
\usepackage[shortlabels]{enumitem}
\newlist{myenumi}{enumerate}{1}
\setlist[myenumi,1]{label=\upshape(\roman*)}
\newlist{myenuma}{enumerate}{1}
\setlist[myenuma,1]{label=\upshape(\alph*)}
\usepackage{color}
\usepackage{todonotes}

\declaretheorem[name=Theorem, numberwithin=section]{thm}
\declaretheorem[name=Lemma,numberlike=thm]{lem}
\declaretheorem[name=Corollary,numberlike=thm]{cor}
\declaretheorem[name=Proposition,numberlike=thm]{prop}
\declaretheorem[name=Definition,numberlike=thm, style=definition]{defi}

\declaretheorem[name=Remark, numberlike=thm, style=remark]{rem}

\numberwithin{equation}{section}
\allowdisplaybreaks[1]
\usepackage[noabbrev]{cleveref} 
\crefname{figure}{Figure}{Figures}
\crefname{table}{Table}{Tables}
\crefname{thm}{Theorem}{Theorems}
\crefname{lem}{Lemma}{Lemmas}
\crefname{defi}{Definition}{Definitions}
\crefname{cor}{Corollary}{Corollaries}
\crefname{prop}{Proposition}{Propositions}
\crefname{conj}{Conjecture}{Conjectures}
\crefname{ex}{Example}{Examples}
\crefname{rem}{Remark}{Remarks}
\crefname{section}{Section}{Sections}
\crefname{chapter}{Chapter}{Chapters}
\crefname{appendix}{Appendix}{Appendices}
\crefdefaultlabelformat{#2\textup{#1}#3}
\creflabelformat{enumi}{(#2#1#3)}
\usepackage{xparse}

\newcommand{\numbers}[1]{\mathbb{#1}}
\newcommand{\C}{\numbers{C}}
\newcommand{\Z}{\numbers{Z}}
\newcommand{\N}{\numbers{N}}
\newcommand{\Q}{\numbers{Q}}
\newcommand{\R}{\numbers{R}}
\newcommand{\ph}{\operatorname{ph}}
\newcommand{\phub}{\underline{\ph}}
\newcommand{\Rm}{\mathcal{R}}
\newcommand{\Fin}{\mathrm{F}}
\newcommand{\Cone}{\operatorname{Cone}}
\newcommand{\SG}{\mathrm{S}}
\newcommand{\ind}{\textup{Ind}}
\newcommand{\Intv}{\textup{I}}

\newcommand{\Sphere}{\mathrm{S}}
\newcommand{\Torus}{\mathrm{T}}

\newcommand{\Cl}{\mathrm{Cl}}

\NewDocumentCommand \Btmfd {} {\mathrm{Bt}}

\NewDocumentCommand \CstarLoc {o} {
  \Cstar_{\mathrm{L}
  \IfValueT{#1}{, #1}
  }
}

\NewDocumentCommand \StolzRel {o d() m m} {
   \mathrm{R}_{
      #3
      \IfValueT{#2}{,(#2)}
    }^{
    \mathrm{spin}
    \IfValueT{#1}{,#1}
   }
   \left(
      #4
   \right)
}

\NewDocumentCommand \StolzPos {o d() m m} {
   \mathrm{Pos}_{
      #3
      \IfValueT{#2}{,(#2)}
    }^{
    \mathrm{spin}
    \IfValueT{#1}{,#1}
   }
   \left(
      #4
   \right)
}

\NewDocumentCommand \StolzPartialPSC {o d() m m} {
   \mathrm{R}\Omega_{
      #3
      \IfValueT{#2}{,(#2)}
    }^{
    \mathrm{spin}
    \IfValueT{#1}{,#1}
   }
   \left(
      #4
   \right)
}

\NewDocumentCommand \SpinBordism {o d() m m} {
   \Omega_{
      #3
      \IfValueT{#2}{,(#2)}
    }^{
    \mathrm{spin}
    \IfValueT{#1}{,#1}
   }
   \left(
      #4
   \right)
}

\NewDocumentCommand \EF {o} {
  \mathrm{E}
  \IfValueT{#1}{_{#1}}
}

\newcommand{\EFin}{\underline{\mathrm{E}}}

\newcommand{\iso}{\cong}

\newcommand{\HZ}{\mathrm{H}}
\newcommand{\KO}{\mathrm{KO}}
\newcommand{\KU}{\mathrm{K}}

\newcommand{\chern}{\operatorname{ch}}
\newcommand{\chernub}{\underline{\chern}}

\newcommand{\sphere}{\mathrm{S}}

\NewDocumentCommand \disk {s} {
  \IfBooleanTF{#1} {
    \mathring{\mathrm{D}}
  } {
    \mathrm{D}
  }
}
\newcommand{\Riem}{\mathcal{R}}

\newcommand*\dif{\mathop{}\!\mathrm{d}}
\newcommand{\tens}{\otimes}
\newcommand{\tensgr}{\mathbin{\widehat{\tens}}}
\newcommand{\Cstar}{\mathrm{C}^*}
\newcommand{\bd}{\partial}
\newcommand{\id}{\operatorname{id}}
\newcommand{\iu}{\mathrm{i}}
\newcommand{\eu}{\mathrm{e}}
\newcommand{\im}{\operatorname{im}}
\newcommand{\OG}{\mathrm{O}}
\newcommand{\SO}{\mathrm{SO}}
\newcommand{\Spin}{\mathrm{Spin}}

\newcommand{\induction}{\operatorname{ind}}
\newcommand{\B}{\mathrm{B}}
\newcommand{\BF}{\B}
\newcommand{\pt}{\ast}

\newcommand{\RSp}{\mathrm{RSp}}
\newcommand{\RU}{\mathrm{RU}}
\newcommand{\RO}{\mathrm{RO}}

\DeclareMathOperator*{\colim}{colim}
\newcommand{\CstarRed}{\Cstar_{\mathrm{r}}}

 \title[The range of the relative higher index and higher rho-invariant]{On the range of the relative higher index and the higher rho-invariant for positive scalar curvature}
 \author{Zhizhang Xie}
 \email{xie@math.tamu.edu}
 \address{Texas A\&M University}

 \author{Guoliang Yu}
 \email{guoliangyu@math.tamu.edu}
 \address{Texas A\&M University}

 \author{Rudolf Zeidler}
 \email{math@rzeidler.eu}
\address{Mathematical Institute, University of Münster, Germany}

\hypersetup{
  colorlinks=true,
  allcolors=blue,
  pdfauthor={Z. Xie, G. Yu, R. Zeidler}
}
\date{}
\begin{document}
 \begin{abstract}
	Let $M$ be a closed spin manifold which supports a positive scalar curvature metric.
	The set of concordance classes of positive scalar curvature metrics on $M$ forms an abelian group $P(M)$ after fixing a positive scalar curvature metric.
	The group $P(M)$ measures the size of the space of positive scalar curvature metrics on $M$. Weinberger and Yu gave a lower bound of the rank of $P(M)$ in terms of the number of torsion elements of $\pi_1(M)$. In this paper, we give a sharper lower bound of the rank of $P(M)$ by studying the image of the relative higher index map from $P(M)$ to the real K-theory of the group $\Cstar$-algebra $\CstarRed(\pi_1(M))$.
	We show that it rationally contains the image of the Baum--Connes assembly map up to a certain homological degree depending on the dimension of $M$.
	At the same time we obtain lower bounds for the positive scalar curvature bordism group by applying the higher rho-invariant.
\end{abstract}
\maketitle
%

\section{Introduction}
Let $M$ be a closed spin manifold. Suppose $M$ carries a positive scalar curvature metric.
The space of all positive scalar curvature metrics on $M$ carries a non-trivial topology, and is in particular usually highly non-connected.
Similarly as in non-existence results for positive scalar curvature metrics, this non-triviality can be detected through the index theory of the spinor Dirac operator. This goes back to the secondary index of Hitchin~\cite{Hitchin74Spinors}.
Index theoretic methods around positive scalar curvature were later enriched by Rosenberg to take the K-theory of the $\Cstar$-algebra of the fundamental group into account~\cite{Rosenberg83CStarPSC}.

A conceptual picture of the interplay between positive scalar curvature and the K-theory of group $\Cstar$-algebras was established by way of mapping the positive scalar curvature bordism sequence of Stolz to the analytic surgery sequence, see \citeauthor{PS14Rho}~\cite{PS14Rho} and \citeauthor{XY14Positive}~\cite{XY14Positive}, and compare \cref{prop:mappingPSCToAnalysis} below.
This involves two types of secondary index invariants --- the \emph{relative higher index} which lies in the K-theory of the group $\Cstar$-algebra and distinguishes connected components in the space of positive scalar curvature metrics, and the \emph{higher rho-invariant} that lies in the analytic structure group and distinguishes bordism classes of positive scalar curvature metrics.
It is a folklore that the image of the higher relative index contains the image of the Baum--Connes assembly map for torsion-free groups, see for instance~\cite{RS01PscSurgery}.
Recently, \citeauthor{ERW17pscFundamentalGroup} obtained such results concerning higher homotopy groups of the space of positive scalar curvature metrics of manifolds with torsion-free fundamental groups~\cite{ERW17pscFundamentalGroup}.
However, no complete results for general groups with torsion have been established so far.
In particular, since the Baum--Connes conjecture predicts that non-trivial higher rho-invariants only exist for groups with torsion, there have been only scarce methods for obtaining positive scalar curvature metrics which can be distinguished up to bordism. Nonetheless, we refer to the work of \citeauthor{BG95Eta}~\cite{BG95Eta} (which deals with finite groups using numerical relative eta-invariants), and \citeauthor{PS07Torsion}~\cite{PS07Torsion} (using the Cheeger-Gromov $L^2$ rho-invariant) for some nice positive results in this direction.

	Following Stolz \cite{Stolz98Concordance}, Weinberger and Yu introduced an abelian group structure, denoted by $P(M)$, on the set of concordance classes of positive scalar curvature metrics on $M$ \cite{WY15FinitePart}. The group $P(M)$ is closely related to the R-group in the exact sequence of Stolz and measures some aspects of the size of the space of positive scalar curvature metrics on $M$.
	Weinberger and Yu then used the finite part of K-theory of the maximal group $\Cstar$-algebra $\Cstar_{\max}(\pi_1(M))$ to give a lower bound of the rank of $P(M)$.\footnote{Note that the statements in~\cite{WY15FinitePart} concerning positive scalar curvature in the dimensions \(4k+1\) are not completely correct if torsion of even order is involved. This is because of an error in a technical proposition~\cite[Proposition 4.4]{WY15FinitePart}. For more details and a remedy of this issue, we refer to~\cref{prop:finCyc,rem:Weinberger-Yu-Error} in the present article.}
	It follows from considerations with the analytic surgery sequence that the elements in $P(M)$ coming from the finite part not only yield different concordance classes but also different positive scalar curvature bordism classes (compare also~\cite{XY17Moduli}).


	More recently, under the assumption that $\pi_1(M)$ satisfies the (rational) strong Novikov conjecture, \citeauthor{BZ17LowDegree}~\cite{BZ17LowDegree} obtained a sharper lower bound of $P(M)$ by incorporating group homology classes of degree up to $2$. In this paper,  under the same assumption that $\pi_1(M)$ satisfies the (rational) strong Novikov conjecture, we shall prove an even sharper lower bound of the rank of $P(M)$ by studying the image of the relative higher index map from $P(M)$ to the K-theory of the reduced group $\Cstar$-algebra $\CstarRed(\pi_1(M))$. Moreover, this new lower bound incorporates \emph{all homology classes up to the dimension of $M$}.
Therefore rationally the entire image of the Baum--Connes assembly map lies in the image of the relative higher index (if we allow the dimension of $M$ to vary).
This extends to full generality what previously has only been known for torsion-free groups.
At the same time, we obtain lower bounds for the size of the image of the higher rho-invariant, and thereby establish a rich source of examples of positive scalar curvature metrics which can be distinguished up to bordism.

The methods in this paper work equally well for the maximal group $\Cstar$-algebra. The maximal version will also give similar applications as the ones stated in the paper. For simplicity, we will only work with the reduced version throughout the paper.

In the following, we will use the universal space for proper actions, denoted by \(\EFin \Gamma\), the universal space for free actions, denoted by \(\EF \Gamma\), and the classifying space \(\BF \Gamma = \EF \Gamma / \Gamma\).
 
Now suppose $M$ is a closed spin manifold of dimension $n \geq 5$ with $\pi_1 M = \Gamma$. Then Stolz' R-group $\StolzRel{n+1}{\BF \Gamma}$ (compare \cref{sec:equivariantStolz} below) acts \emph{freely and transitively} on the set of concordance classes of positive scalar curvature metrics on $M$, which is denoted by $ \tilde{\pi}_0(\Riem^+(M))$, see \cite[Theorem~1.1]{Stolz98Concordance}.
After choosing a base point $g_0$ of $ \tilde{\pi}_0(\Riem^+(M))$, there is a bijection between  $\StolzRel{n+1}{\BF \Gamma}$  and $ \tilde{\pi}_0(\Riem^+(M))$. In particular, this bijection introduces an abelian group structure on $ \tilde{\pi}_0(\Riem^+(M))$.
In fact, it is not difficult to see that $\StolzRel{n+1}{\BF \Gamma}$ is isomorphic to the group $P(M)$ of \citeauthor{WY15FinitePart} (see \cref{sec:twoDefsOfConcordanceGp}).

The relative higher index map
\[ \alpha\colon \StolzRel{n+1}{\BF \Gamma} = P(M) \to \KO_{n+1}(\CstarRed(\Gamma))\] takes a positive scalar curvature metric $g$ to the relative higher index $\ind_\Gamma(g, g_0)$, that is, the higher index of the Dirac operator on the cylinder $M\times \R$, where $M \times \R$ is equipped with  a Riemannian metric $g_t + {\dif t}^2$ such that $g_t = g_0$ for $t \leq 0$ and $g_t = g$ for $t \geq 1$.
It follows from the relative higher index theorem  that $\alpha$ is a group homomorphism.
A lower bound for the rank of the image of $\alpha$ will also serve as a lower bound of the rank of $P(M)$.

The higher rho-invariant is a homomorphism
\[\rho \colon \StolzPos{n}{\BF \Gamma} \to \SG_{n}^\Gamma(\EF \Gamma)\]
from the positive scalar curvature bordism group (which appears in Stolz' sequence) to the \emph{analytic structure group} (which appears in the analytic surgery sequence).
Using the comparison diagram between Stolz' sequence and the analytic surgery sequence (see \cref{prop:mappingPSCToAnalysis}), one can deduce lower bounds for the image of \(\rho\) from lower bounds for the image of \(\alpha\).

In order to estimate the size of the image of the relative higher index map $\alpha$, we introduce a new bordism group $\StolzPartialPSC[\Gamma]{n+1}{E_+, E_-}$ for each pair of proper $\Gamma$-spaces $E_-\subseteq E_+$. Roughly speaking, it is a hybrid of the Stolz' group  $\StolzRel[\Gamma]{n+1}{E_-}$ and the equivariant spin bordism group $\SpinBordism[\Gamma]{n+1}{E_+}$.
It admits a natural higher index map $\ind_\Gamma \colon \StolzPartialPSC[\Gamma]{n}{E_+, E_-} \to \KO_{n}(\CstarRed \Gamma)$, which factors through the relative higher index map $\alpha \colon \StolzRel[\Gamma]{n+1}{E_-}\to \KO_{n+1}(\CstarRed(\Gamma))$, see \cref{cor:factorsThrough} below. In particular, it follows that the image of $\alpha$ contains the image of the map $\ind_\Gamma \colon \StolzPartialPSC[\Gamma]{n+1}{\EFin \Gamma, \EF \Gamma} \to \KO_{n+1}(\CstarRed(\Gamma))$.
As a consequence, we are reduced to studying the image of $\ind_\Gamma \colon \StolzPartialPSC[\Gamma]{n+1}{\EFin \Gamma, \EF \Gamma} \to \KO_{n+1}(\CstarRed(\Gamma))$.
One main theorem of the paper is to give a lower estimate of the image of this map.
To this end, we use the the \emph{equivariant delocalized Pontryagin character} to identify the K-homology group $\KO_{p}^\Gamma(\EFin \Gamma) \tens \C$ with the following expression in terms of ordinary group homology,
\[
  \bigoplus_{k \in \Z} \HZ_{p+4k}(\Gamma; \Fin^0\Gamma) \oplus \HZ_{p - 2 + 4k}(\Gamma; \Fin^1 \Gamma),
\] 
where \(\Fin^0 \Gamma\) and \(\Fin^1 \Gamma\) denotes the symmetric and anti-symmetric part, respectively, of the vector space generated by finite order elements of \(\Gamma\); see~\cref{subsec:chern} for details on this construction.
\begin{thm}
	Let $p \in \{0,1,2,3\}$ and $k \geq 1$.
	Then the image of the map
	\begin{equation*}
	\ind_\Gamma\colon \StolzPartialPSC[\Gamma]{4k+p}{\EFin \Gamma, \EF \Gamma} \tens \C  \to  \KO_{p}(\CstarRed \Gamma) \tens \C
	\end{equation*}
	contains
	\begin{equation*}
	\mu \left( \bigoplus_{l < k} \HZ_{4l + p}(\Gamma; \Fin^0\Gamma) \oplus \HZ_{4l - 2 + p}(\Gamma; \Fin^1 \Gamma) \right)
	\end{equation*}
\end{thm}
Here $\mu \colon \KO_{p}^\Gamma(\EFin \Gamma) \to \KO_p(\CstarRed \Gamma)$ denotes the real version of the assembly map that features in the Baum--Connes conjecture~\cite{BC00Conjecture}.
Note that the complex version of the conjecture implies the real version~\cite{BK04Real}.
After inverting \(2\) (in particular, rationally), injectivity and surjectivity of the complex Baum--Connes assembly map is separately equivalent to the corresponding statement for the real Baum--Connes assembly map~\cite{Schick:RealVsComplex}.
One key ingredient for the proof of the theorem above is the realization of (rational) $\KO$-homology classes by spin $\Gamma$-manifolds with appropriate control over their dimensions.

As an application, our new lower bound for the rank of $\StolzRel{n+1}{\BF \Gamma} = P(M)$ follows  immediately from the theorem above.
\begin{cor}
	Let $p \in \{0,1,2,3\}$ and $k \geq 1$.
	Suppose that the Baum--Connes assembly map $\mu \colon \KO_\ast^\Gamma(\EFin \Gamma) \to \KO_\ast(\CstarRed \Gamma)$ is rationally injective.
	Then the rank of $\StolzRel{4k + p}{\BF \Gamma}$ is at least the dimension of
	\begin{equation*}
	\bigoplus_{l < k} \bigoplus_{q \in \{0,1\}} \HZ_{4l -2q + p}(\Gamma;  \Fin^q\Gamma).
	\end{equation*}
  In particular, if  $\dim M = 4k+p - 1 \geq 5$, then the same lower bound applies to \(P(M)\).
\end{cor}

In the presence of upper bounds on the rational homological dimension and surjectivity of the Baum--Connes assembly map, our method yields surjectivity of the relative index map:

\begin{cor}
  Let $n \geq 4$ and the rational homological dimension of $\Gamma$ be at most $n -3$. Suppose that the Baum--Connes assembly map for $\Gamma$ is rationally surjective.
  Then the relative index map $\alpha \colon \StolzRel{n}{\BF \Gamma} \to \KO_{n}(\CstarRed \Gamma)$ is rationally surjective.
\end{cor}

Similar conclusions apply to the positive scalar curvature bordism group:
\begin{cor}
      Let $p \in \{0,1,2,3\}$ and $k \geq 1$.
    Suppose that the Baum--Connes assembly map $\mu \colon \KO_\ast^\Gamma(\EFin \Gamma) \to \KO_\ast(\CstarRed \Gamma)$ is rationally injective.
    Then the rank of $\StolzPos{4k + p - 1}{\BF \Gamma}$ is at least the dimension of
    \begin{equation*}
      \bigoplus_{l < k} \bigoplus_{q \in \{0,1\}} \HZ_{4l -2q + p}(\Gamma;  \Fin^q_0\Gamma).
    \end{equation*}
\end{cor}

\begin{cor}
  Let $n \geq 4$ and the rational homological dimension of $\Gamma$ be at most $n -3$. Suppose that the Baum--Connes assembly map for $\Gamma$ is a rational isomorphism.
  Then the higher rho-invariant \(\rho \colon \StolzPos{n-1}{\BF \Gamma} \to \SG_{n-1}^\Gamma(\EF \Gamma)\) is rationally surjective.
\end{cor}

Note that even for finite groups and groups of rational homological dimension at most \(2\), our results above are stronger than what can be obtained from \citeauthor{BG95Eta}~\cite{BG95Eta} and \citeauthor{BZ17LowDegree}~\cite{BZ17LowDegree} because our estimates begin in dimension \(4\), whereas these previous results only apply to dimension \(6\) and above.

The principal reason why assumptions on homological dimension are necessary in order to obtain full surjectivity results is that KO-theory is $8$-periodic, whereas the various bordism groups which appear in Stolz' sequence are not.
However, this can be remedied by force by introducing Bott periodicity formally.
Indeed, let $\Btmfd$ denote the \emph{Bott manifold}, that is, an $8$-dimensional simply connected spin manifold with A-hat genus $\hat{\mathrm{A}}(\Btmfd) =1$.
Then let
\begin{equation}
  \StolzRel{n}{\BF \Gamma}\left[\Btmfd^{-1}\right] \coloneqq \colim \left(\StolzRel{n}{\BF \Gamma} \xrightarrow{\times \Btmfd} \StolzRel{n+8}{\BF \Gamma} \xrightarrow{\times \Btmfd} \dotsm \right). \label{eq:bottStabilize}
  \end{equation}
By Bott periodicity the relative index map induces the stabilized map
\begin{equation}
  \alpha[\Btmfd^{-1}] \colon \StolzRel{n}{\BF \Gamma}\left[\Btmfd^{-1}\right] \to \KO_n(\CstarRed \Gamma).\label{eq:stabilizedRelIndex}
\end{equation}
Similarly, there is a Bott-stabilized version of the higher rho-invariant
\begin{equation}
  \rho[\Btmfd^{-1}] \colon \StolzPos{n-1}{\BF \Gamma}[\Btmfd^{-1}] \to \SG_{n-1}(\EF \Gamma).\label{eq:stabilizedRho}
\end{equation}
Our main results imply:
\begin{cor} \label{cor:stabilized}
   Suppose that the Baum--Connes assembly map for $\Gamma$ is rationally surjective.
  Then the stabilized relative index map \labelcref{eq:stabilizedRelIndex} is rationally surjective.
  If the Baum--Connes assembly map for \(\Gamma\) is a rational isomorphism, then the stabilized rho-invariant \labelcref{eq:stabilizedRho} is also rationally surjective.
\end{cor}

However, note that while this stabilization procedure conceptually suggests itself and is necessary to get a statement as in \cref{cor:stabilized} via our methods, it is not clear that it is required in principle.
Indeed, unlike in the case of the primary index, there is no known obstruction that would preclude the relative higher index of two positive scalar curvature metrics from being a K-theory class associated to a K-homology class of a homological degree higher than the underlying manifold.
Neither are there any examples of this kind.
Moreover, in case the Baum--Connes conjecture fails, there is also no a priori reason why the relative index needs to be in the image of the assembly map.
This means that, while our present results rationally exhaust what is possible through known geometric constructions together with the Baum--Connes conjecture, it remains a tantalizing open question whether the space of secondary index invariants associated to positive scalar curvature contains any more exotic elements. 

The paper is organized as follows. In \cref{sec:equivariantStolz}, we introduce the hybrid bordism group of Stolz' R-group and the spin bordism group, and show that its higher index map factors through the relative higher index map. In \cref{sec:mainResults}, we show how to realize (rational) $\KO$-homology classes by spin $\Gamma$-manifolds with appropriate control over their dimensions, then apply it to study the image of the higher index map on the hybrid bordism group introduced in \cref{sec:equivariantStolz}.
We then apply these results to obtain a sharper lower bound of the rank of Stolz' R-group and the positive scalar curvature bordism group.
In \cref{sec:twoDefsOfConcordanceGp}, we show that the two definitions of the group of concordance classes of positive scalar curvature metrics agree.

\section{The equivariant positive scalar curvature sequence}\label{sec:equivariantStolz}
In this section, we review the equivariant version of the positive scalar curvature sequence of Stolz and introduce a new relative group that interpolates between the equivariant spin bordism group and the non-equivariant version of Stolz' $\mathrm{R}$-group.
This group is the main new conceptual tool developed in the present paper and will be used in the proof of our main results in \cref{sec:mainResults}.

\begin{defi}[{Compare \cite[Section 5]{XY14Positive}}]
  Let $\Gamma$ be a discrete group and $E$ a proper $\Gamma$-space.
  Then we have an equivariant version of Stolz' positive scalar curvature seqence
  \begin{equation}
  \resizebox{0.9\linewidth}{!}{$
    \cdots \rightarrow
    \StolzRel[\Gamma]{n+1}{E} \xrightarrow{\bd} \StolzPos[\Gamma]{n}{E}
      \xrightarrow{q} \SpinBordism[\Gamma]{n}{E} \xrightarrow{j} \StolzRel[\Gamma]{n}{E}
      \rightarrow \cdots,$}
  \end{equation}
  which is defined in the same way as the standard sequence of Stolz~\cite{Stolz98Concordance} except for replacing compact spin manifolds with proper $\Gamma$-cocompact spin manifolds, and continuous maps by $\Gamma$-equivariant continuous maps, everywhere.

  More precisely, $\SpinBordism[\Gamma]{n}{E}$ consists of $\Gamma$-equivariant spin bordism classes of pairs $(M, \phi)$, where $M$ is a proper cocompact spin $\Gamma$-manifold without boundary and $\phi \colon M \to E$ an equivariant continuous map.
  The group $\StolzPos[\Gamma]{n}{E}$ consists of equivariant spin bordism classes of triples $(M, \phi, g)$, where $(M, \phi)$ is as before and $g \in \Riem^+(M)^\Gamma$ is a $\Gamma$-invariant metric of uniformly positive scalar curvature.
  Here bordisms are required to be endowed with $\Gamma$-invariant metrics of uniformly positive scalar curvature which are collared near the boundaries.
  Finally, $\StolzRel[\Gamma]{n+1}{E}$ consists of suitable equivariant bordism classes of triples $(W, \phi, g)$, where $W$ is a proper cocompact spin $\Gamma$-manifold with boundary, $\phi \colon W \to E$ an equivariant continuous map, and $g \in \Riem^+(\bd W)^\Gamma$ a $\Gamma$-invariant metric of uniformly positive scalar curvature at the boundary.
  The maps $\bd$, $q$ and $j$ are the evident forgetful maps.
\end{defi}

\begin{prop}[{\cite[Theorem B]{XY14Positive}, see also~\cite{PS14Rho,Z16Positive,Z16PhD}}]\label{prop:mappingPSCToAnalysis}
  There is a map of exact sequences mapping the equivariant positive scalar curvature sequence to analysis:
  \begin{equation}
  \resizebox{0.9\linewidth}{!}{
     \begin{tikzcd}[column sep=small, ampersand replacement=\&]
        \SpinBordism[\Gamma]{n}{E} \rar \dar{\beta} \&
        \StolzRel[\Gamma]{n}{E} \rar{\partial} \dar{\alpha} \&
        \StolzPos[\Gamma]{n-1}{E} \rar \dar{\rho} \&
        \SpinBordism[\Gamma]{n-1}{E} \rar \dar{\beta} \&
        \StolzRel[\Gamma]{n-1}{E} \dar{\alpha}
        \\
        \KO_{n}^\Gamma(E) \rar{\mu} \&
        \KO_{n}(\CstarRed \Gamma) \rar{\partial}
          \& \SG_{n-1}^\Gamma(E) \rar
          \& \KO_{n-1}^\Gamma(E) \rar{\mu}
          \& \KO_{n-1}(\CstarRed \Gamma)
     \end{tikzcd}}\label{eq:mappingPSCtoAnalysis}
  \end{equation}
\end{prop}

Here the bottom horizontal sequence is the real version of the analytic surgery sequence of Higson and Roe~\cite{HR05MappingI,HR05MappingII,HR05MappingIII}, where $\SG_{\ast}^\Gamma(E)$ denotes the \emph{analytic structure group}.

The relevant case in the study of positive scalar curvature on closed manifolds with fundamental group $\Gamma$ is $E = \EF \Gamma$.
Therefore we are interested in studying the size of the image of $\alpha$ and $\rho$ for $E = \EF \Gamma$.
However, the Baum--Connes conjecture predicts that in order to understand $\KO_\ast(\CstarRed \Gamma)$, we need to consider $E = \EFin \Gamma$ instead.
Thus our goal is to develop a method for lifting data from $\KO_\ast^\Gamma(\EFin \Gamma)$ to $\StolzRel[\Gamma]{\ast}{\EF \Gamma}$.
To that end, we construct a new relative bordism group $\StolzPartialPSC[\Gamma]{n}{\EFin \Gamma, \EF \Gamma}$ à la Stolz which comes with a natural map to both $\StolzRel[\Gamma]{\ast}{\EF \Gamma}$ and $\SpinBordism[\Gamma]{n}{\EFin \Gamma}$.
Then our strategy in the next section is to use the equivariant Chern character to show that data from $\KO_\ast^\Gamma(\EFin \Gamma)$ can be lifted to $\StolzPartialPSC[\Gamma]{n}{\EFin \Gamma, \EF \Gamma}$.

Roughly speaking, $\StolzPartialPSC[\Gamma]{n}{\EFin \Gamma, \EF \Gamma}$ consists of proper $\Gamma$-spin manifolds which are partitioned by a codimension $1$ hypersurface together with a positive scalar curvature metric on one half and such that the group action is free on the other half.
The following defintion makes this precise.

\begin{defi}\label{defi:StolzPartialPSC}
  Let $\Gamma$ be a discrete group and $E_- \subseteq E_+$ be a pair of proper $\Gamma$-spaces.
  We define a group $\StolzPartialPSC[\Gamma]{n}{E_+, E_-}$ of bordism classes as follows:
  \begin{itemize}
    \item   Its elements are represented by tuples of the form $(W_-, W_+, \sigma_\mp, g)$, where $W_\mp$ are complete spin $n$-manifolds endowed with a cocompact $\Gamma$-action and $\bd W_+ = \bd W_- =: M$.
      Moreover, $\sigma_\mp \colon W_\mp \to E_\mp$ are continuous $\Gamma$-equivariant maps which agree on $M$, and $g \in \Rm^+(W_+)$ is a $\Gamma$-invariant metric of positive scalar curvature (collared near $\bd W_+$).
    \item The tuples $(W_-, W_+, \sigma_\mp, g)$ and $(W_-^\prime, W_+^\prime, \sigma^\prime_\mp, g^\prime)$ are \emph{bordant} if the following holds:
    \begin{itemize}
      \item There exists spin bordisms $V_\mp$ between $W_\mp$ and $W_\mp^\prime$ which restrict to the same bordism between $M = \bd W_+ = \bd W_-$ and $M^\prime = \bd W_+^\prime = \bd W_-^\prime$.
      More precisely, $V_\mp$ is a spin manifold with corners and its boundary decomposes as $\bd V_\mp = W_\mp \cup_M N \cup_{M^\prime} W_\mp^\prime$, where $N$ is a spin bordism between $M$ and $M^\prime$
      \item There exist maps $S \colon V_\mp \to E_\mp$ which restrict to $\sigma_\mp \sqcup \sigma^\prime_\mp$ and agree on $N$.
      \item There exists a metric $h \in \Rm^+(V_+)$ which is collared near the boundary and restricts to $g \sqcup g^\prime$ on $W_+ \sqcup W_+^\prime$.
    \end{itemize}
  \end{itemize}

\end{defi}

\begin{defi}\label{rem:partialPSCForget}
  We define the following forgetful maps:
  \begin{align*}
    \StolzPartialPSC[\Gamma]{n}{E_+, E_-} &\xrightarrow{r} \StolzRel[\Gamma]{n}{E_-},\
    (W_-, W_+, \sigma_\mp, g) \mapsto (W_-, \sigma_-, g|_{\bd W_+}) \\
    \StolzPartialPSC[\Gamma]{n}{E_+, E_-} &\xrightarrow{\omega} \SpinBordism[\Gamma]{n}{E_+},\
    (W_-, W_+, \sigma_\mp, g) \mapsto (W_- \cup_M W_+, \sigma_- \cup \sigma_+)
  \end{align*}
\end{defi}

\begin{rem}[Induction] \label{rem:induction}
  Let $E_- \subseteq E_+$ be a pair of $G$-spaces.
  Let $\psi \colon G \to \Gamma$ be a group homomorphism such that $\ker(\psi)$ acts freely on $E_+$.
  Then there is an induction map
  \begin{equation*}
    \induction_\psi \colon \StolzPartialPSC[G]{n}{E_+, E_-} \to \StolzPartialPSC[\Gamma]{n}{\Gamma \times_\psi E_+, \Gamma \times_\psi E_-}
  \end{equation*}
  taking a cycle represented by $(W_-, W_+, \sigma_\mp, g)$ to the cycle represented by $(\Gamma \times_\psi W_-, \Gamma \times_\psi W_+,  \id_\Gamma \times_\psi \sigma_\mp, \Gamma \times_\psi g)$.
\end{rem}

\begin{prop}\label{prop:interpolatingDiagram}
  Let $E_- \subseteq E_+$ be a pair of proper $\Gamma$-spaces.
  Then we have a commutative diagram:
  \begin{equation*}
    \begin{tikzcd}
      \StolzPartialPSC[\Gamma]{n}{E_+, E_-} \rar["r"] \dar["\omega"] & \StolzRel[\Gamma]{n}{E_-} \dar["\iota_\ast"] \\
      \SpinBordism[\Gamma]{n}{E_+} \rar["i"] & \StolzRel[\Gamma]{n}{E_+}
    \end{tikzcd}
  \end{equation*}
\end{prop}
\begin{proof}
\begin{sloppypar}
  Let $x \in \StolzPartialPSC[\Gamma]{n}{E_+, E_-}$.
  Then $x$ is represented by a tuple $(W_-, W_+, \sigma_\pm, g)$ as in \cref{defi:StolzPartialPSC}.
  Let $W := W_- \cup_M W_+$ denote the proper $\Gamma$-spin manifold obtained by gluing together $W_-$ with $W_+$ along $M := \bd W_- = \bd W_+$.
  This manifold admits a map $\sigma := (\iota \circ \sigma_-) \cup \sigma_+ \colon W \to E_+$, where $\iota \colon E_- \hookrightarrow E_+$ denotes the inclusion map.
  Then $i(\omega(x)) \in \StolzRel[\Gamma]{n}{E_+}$ is represented by the tuple $(W, \sigma, \emptyset)$, where $\emptyset$ stands for the \enquote{empty metric} ($W$ does not have a boundary).
  On the other hand, the element $\iota_\ast(r(x))$ is represented by $(W_-, \iota \circ \sigma_-, g|_{M})$.
  \end{sloppypar}

 Now we can view $W \times [0,1]$ as a bordism between $W$ and $W_-$. To see this, we interpret one copy of $W_+$ inside $\bd(W \times [0,1])$ as a bordism between $\bd W = \emptyset$ and $\bd W_- = \bd W_+ = M$, and write $\bd(W \times [0,1])$ as the union $W \sqcup W_+ \cup_{M} W_-$ (compare \cref{fig:bordism}).
  By construction, the positive scalar curvature metric $g|_{M}$ on $\bd W_-$ extends to the $\Gamma$-invariant positive scalar curvature metric $g$ on $W_+$.
  Moreover, the map $(\iota \circ \sigma_-) \sqcup \sigma$ on $W_- \sqcup W$ extends to a map on $W \times [0,1]$ via $(p,t) \mapsto \sigma(p)$.
  Thus we have constructed a bordism between $(W, \sigma, \emptyset)$ and $(W_-, \iota \circ \sigma_-, g|_{M})$ which witnesses that they represent the same element of $\StolzRel[\Gamma]{n}{E_+}$.
  This proves $i(\omega(x)) = \iota_\ast(r(x))$.
  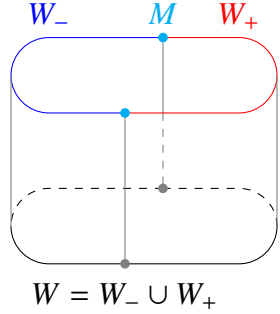
\begin{figure}[h]
 \caption{The bordism $W \times [0,1]$.}\label{fig:bordism}
  \begin{tikzpicture}
  \coordinate (M0) at (0cm, 0cm);
  \coordinate (M2) at ($(M0) + (0cm, 2cm)$);
  \coordinate (M1) at (0.5cm, 1cm);
  \coordinate (M3) at ($(M1) + (0cm, 2cm)$);

  \coordinate (LB0) at (-1cm, 0cm);
  \coordinate (LB1) at (-1.5cm, 0.5cm);
  \coordinate (LB2) at (-1cm, 1cm);

  \coordinate (RB0) at (1.5cm, 0cm);
  \coordinate (RB1) at (2cm, 0.5cm);
  \coordinate (RB2) at (1.5cm, 1cm);

  \coordinate (LT0) at ($(LB0) + (0cm, 2cm)$);
  \coordinate (LT1) at ($(LB1) + (0cm, 2cm)$);
  \coordinate (LT2) at ($(LB2) + (0cm, 2cm)$);

  \coordinate (RT0) at ($(RB0) + (0cm, 2cm)$);
  \coordinate (RT1) at ($(RB1) + (0cm, 2cm)$);
  \coordinate (RT2) at ($(RB2) + (0cm, 2cm)$);

  \coordinate (RTI) at
	  (intersection cs: first line={(M2) -- (RT0)},
			    second line={(M3) -- (M1)});

	\draw [gray] (M3) -- (RTI);
	\draw [gray, dashed] (RTI) -- (M1);

	\draw [dashed] (M1) -- (LB2) arc (90:180:0.5cm);
	\draw (LB1) arc (180:270:0.5cm) -- (M0);
	\draw [blue] (M3) -- (LT2) arc (90:270:0.5cm) -- (M2);

	\draw (M0) -- (RB0) arc (-90:0:0.5cm);
	\draw [dashed] (RB1) arc (0:90:0.5cm) -- (M1);
	\draw [red] (M2) -- (RT0) arc (-90:90:0.5cm) -- (M3);

	\draw [gray] (LB1) -- (LT1);
	\draw [gray] (RB1) -- (RT1);

	\draw [gray] (M2) -- (M0);

	  \foreach \i in {2,3}
  {
	  \fill [cyan] (M\i) circle (0.15em);
	}

		  \foreach \i in {0,1}
  {
	  \fill [gray] (M\i) circle (0.15em);
	}

	\node [above=0.6cm, blue, anchor=north] at (LT2) {$W_-$};
	\node [above=0.6cm, cyan, anchor=north] at (M3) {$M$};
	\node [above=0.6cm, red, anchor=north] at (RT2) {$W_+$};
	\node [below=0.75cm, black, anchor=south] at (M0) {$W = W_- \cup W_+$};
\end{tikzpicture}
\end{figure}
  \end{proof}

\begin{rem}
  We will apply this construction to $E_- = \EF \Gamma$ and $E_+ = \EFin \Gamma$, where we ensure that $\EF \Gamma$ is a subspace of $\EFin \Gamma$ by passing to a mapping cylinder if necessary.
  In this case, $\StolzRel[\Gamma]{n}{\EF \Gamma} = \StolzRel{n}{\BF \Gamma}$.
\end{rem}

\begin{cor}\label{cor:factorsThrough}
  The image of the relative index map $\alpha \colon \StolzRel{n}{\BF \Gamma} \to \KO_n(\CstarRed \Gamma)$ contains the image of the composition
  \begin{equation*}
    \ind_\Gamma \colon \StolzPartialPSC[\Gamma]{n}{\EFin \Gamma, \EF \Gamma} \xrightarrow{\omega} \SpinBordism[\Gamma]{n}{\EFin \Gamma} \xrightarrow{\beta} \KO_n^\Gamma(\EFin \Gamma)  \xrightarrow{\mu} \KO_n(\CstarRed \Gamma).
  \end{equation*}
\end{cor}
\begin{proof}
  As a consequence of \cref{prop:interpolatingDiagram} we obtain the commutative diagram
    \begin{equation*}
    \begin{tikzcd}
      \StolzPartialPSC[\Gamma]{n}{\EFin \Gamma, \EF \Gamma} \rar["r"] \dar["\omega"] & \StolzRel[\Gamma]{n}{\EF \Gamma} \dar["\iota_\ast"]  \rar[equal] & \StolzRel{n}{\BF \Gamma} \ar[ddl, bend left, "\alpha"] \\
      \SpinBordism[\Gamma]{n}{\EFin \Gamma} \rar["i"] \dar["\beta"]
        & \StolzRel[\Gamma]{n}{\EFin \Gamma} \dar["\alpha"] \\
      \KO_n^\Gamma(\EFin \Gamma) \rar["\mu"] & \KO_n(\CstarRed \Gamma).
    \end{tikzcd}
  \end{equation*}
  This shows that the composition of maps in the statement of the corollary factors through the relative index map $\alpha \colon \StolzRel{n}{\BF \Gamma} \to \KO_n(\CstarRed \Gamma)$.
\end{proof}

\section{The image of the relative higher index maps}\label{sec:mainResults}
In this section, we prove our main results.
We start by reviewing relevant material on the equivariant delocalized Chern character in~\cref{subsec:chern}.
Then we collect statements for geometrically realizing (equivariant) KO-homology classes in~\cref{subsec:geom}.
Finally, we put these ingredients together to lift equivariant KO-homology classes to the group $\StolzPartialPSC[\Gamma]{\ast}{\EFin \Gamma, \EF \Gamma}$ which was defined in \cref{sec:equivariantStolz}.
From this we deduce new quantitative lower bounds on the rank of Stolz' groups $\StolzRel{\ast}{\BF \Gamma}$ and $\StolzPos{\ast}{\BF \Gamma}$.
In particular, we thereby substantially extend earlier results of \citeauthor{BG95Eta}~\cite{BG95Eta} and \citeauthor{BZ17LowDegree}~\cite{BZ17LowDegree}.
\subsection{Chern characters}\label{subsec:chern}
If $X$ is a space, there is the homological Chern character $\chern \colon \KU_p(X) \tens \Q \xrightarrow{\cong} \bigoplus_{k \in \Z} \HZ_{p + 2k}(X; \Q)$ for $p \in \{0,1\}$.
Precomposing this with complexification yields the Pontryagin character $\ph \colon \KO_p(X) \tens \Q \xrightarrow{\cong} \bigoplus_{k \in \Z} \HZ_{p + 4k}(X; \Q)$ for $p \in \{0,1,2,3\}$ (using that real K-homology is rationally $4$-periodic).
In particular, if $\Gamma$ is a group, we obtain the Pontryagin character for group homology $\ph \colon \KO_p^\Gamma(\EF \Gamma) \tens \Q \iso \KO_p(\B \Gamma) \tens \Q \xrightarrow{\cong} \bigoplus_{k \in \Z} \HZ_{p + 4k}(\Gamma; \Q)$.

Next we turn to the equivariant setting for proper actions.
Let $\Fin \Gamma$ be the complex vector space generated freely by all the finite order elements in $\Gamma$.
The action by conjugation of $\Gamma$ on $\Fin \Gamma$ turns $\Fin \Gamma$ into a $\C \Gamma$-module.
The \emph{equivariant delocalized Chern character}, first introduced by \citeauthor{BC88Chern}~\cite{BC88Chern}, yields an isomorphism
\begin{equation*}
  \chernub_\Gamma \colon \KU_p^\Gamma(\EFin \Gamma) \tens \C \xrightarrow{\cong}
  \bigoplus_{k \in \Z}\HZ_{p +2k}(\Gamma; \Fin \Gamma), \qquad p \in \Z.
\end{equation*}
For technical purposes, we will work with the \enquote{handicrafted Chern character} of Matthey~\cite[Theorem~1.4]{matthey:delocChern}.
Moreover, to use it for our applications, we need a real version of it.
We have that $\Fin \Gamma = \Fin^0 \Gamma \oplus \Fin^1 \Gamma$, where $\Fin^q \Gamma = \{ f \in \Fin \Gamma \mid f(\gamma^{-1}) = (-1)^q f(\gamma) \}$.
Precomposing the delocalized equivariant Chern character with the complexification map $\KO_p^\Gamma(\EFin \Gamma) \to \KU_p^\Gamma(\EFin \Gamma)$ yields the \emph{equivariant delocalized Pontryagin character} (see \cite[Section 2]{BZ17LowDegree} for more details):
\begin{equation*}
\resizebox{\linewidth}{!}{$
\phub_\Gamma \colon \KO_p^\Gamma(\EFin \Gamma) \tens \C \xrightarrow{\cong}
\bigoplus_{k \in \Z} \HZ_{p+4k}(\Gamma; \Fin^0\Gamma) \oplus \HZ_{p - 2 + 4k}(\Gamma; \Fin^1 \Gamma),\ p \in \Z.$}
\end{equation*}

\subsection{Geometric ingredients}\label{subsec:geom}
Now we state the main ingredients for geometrically realizing $\KO$-homology classes.
The first is a folklore concerning rational homology:

\begin{prop}\label{prop:framedBordism}
  Let $X$ be a space and $n \geq 0$.
  Then $\HZ_n(X; \Q)$ is generated by Pontryagin characters of the spin fundamental classes of closed $n$-dimensional spin manifolds.
\end{prop}
\begin{proof}
  Consider framed bordism, that is, the bordism homology theory of stably
parallelizable manifolds. The framed bordism groups of a point are
isomorphic to the stable homotopy groups of spheres by the original Pontryagin--Thom isomorphism.
In particular, the framed bordism homology rationally agrees with ordinary homology.
So every rational homology class can be realized as a rational multiple of a class represented by a
(stably) parallelizable manifold, which is in particular spin.
Moreover, if $M$ is a (stably) parallelizable manifold, then its Pontryagin classes must vanish except in degree $0$.
Hence the Pontryagin character of its spin fundamental class is the same as its homological fundamental class.
\end{proof}

The next theorem implies that rationally all information in equivariant KO-homology coming from the representation theory of finite cyclic groups can be realized by $4$- and $6$-dimensional equivariant spin manifolds.

\begin{thm}\label{prop:finCyc}
  Let $m \geq 2$, $k \geq 1$ and $q \in \{0,1\}$.
  The group $\KO_{2q}^{\Z/m}(\pt) \tens \Q$ is generated by equivariant indices of compact spin $(\Z/m)$-manifolds of dimension $4k+2q$, where the action is free outside an invariant submanifold of codimension $4$ with trivial normal bundle.
\end{thm}
\begin{rem}\label{rem:Weinberger-Yu-Error}
A similar claim was made by \citeauthor{WY15FinitePart} in \cite[Proposition~4.4]{WY15FinitePart}.
In the case \(q = 0\), the statement in loc.~cit.\ essentially agrees with the formulation given here.
However, in the case \(q = 1\) and \(m\) an even integer, the original statement of \citeauthor{WY15FinitePart} is not correct.
The reason for this is---as we shall discuss in the proof below---that after inverting the prime \(2\), the group $\KO_{2}^{\Z/m}(\pt)$ can be identified with the group generated by antisymmetric characters on the cyclic group \(\Z/m\).
For \(m\) even, the rank of this group is actually one less than what is claimed in~\cite[Proposition~4.4]{WY15FinitePart}.
This can be seen most acutely in the case of a cyclic group of order two because \(\KO_2^{\Z/2}(\ast) \otimes \Q = 0\).
Indeed, the example given in the proof of \cite[Proposition~4.4]{WY15FinitePart} does not work for the base case \(m = 2\) and so the inductive argument cannot get started for even \(m\).
The arguments given for the other cases (\(q = 0\) or \(m\) odd) conceivably still work and would give the correct result but the full details of the equivariant index computation were not given.
To remedy all of this, we provide a new proof here. It is partly inspired by the construction in \cite{BG95Eta}.
\end{rem}
\begin{proof}[Proof of~\cref{prop:finCyc}]
  In the course of the proof, we will have to relate equivariant real K-theory to the representation theory of finite groups.
  In the following, we will give a brief account of the necessary identifications and computations; for more details we refer to~\cite[Chapter~2]{BG10Real}.
  To this end, let \(H\) be a finite group and write \(\RO(H)\), \(\RU(H)\) and \(\RSp(H)\) for the Grothendieck groups of finite dimensional real, complex and quaternionic representations, respectively.
  There are complexification maps \(\RO(H) \to \RU(H)\) and \(\RSp(H) \to \RU(H)\).
  In these terms, there are canonical identifications $\KO_0^{H}(\pt) = \RO(H)$ and $\KO_2^{H}(\pt) = \RU(H)  / \RSp(H)$, see for instance \cite[Theorem~2.2.12]{BG10Real}.
  After inverting $2$ the latter is the same as $\RU(H)/(1+\tau)$, where $\tau \colon \RU(H) \to \RU(H)$ is the map induced by complex conjugation.
  Moreover, after inverting $2$ the group $\RU(H)/(1+\tau)$ can be identified with the image of $(1-\tau)$ in $\RU(H)$.
  Similarly, $\RO(H)$ can be identified with the image of $(1+\tau)$.
  We thus obtain a decomposition
  \[ \RU(H)\left[\frac{1}{2}\right] = \im(1 + \tau) \oplus \im(1 - \tau),\]
  where the first summand corresponds to  $\KO^H_0(\ast)\left[\frac{1}{2}\right]$ and the second to $\KO^H_2(\ast)\left[\frac{1}{2}\right]$.
  Note that $\im (1+\tau)$ is the subgroup generated by those virtual representations with symmetric characters and $\im (1-\tau)$ is generated by those with antisymmetric characters.

 It suffices to prove the proposition for \(k=1\), that is, produce \(4\)- and \(6\)-dimensional manifolds.
  To obtain higher-dimensional examples, we can then simply take products with copies of the particular Kummer surface \(V = \{[Z_0 : Z_1 : Z_2 : Z_3] \in \C\mathrm{P}^3 \mid Z_0^4 + Z_1^4 + Z_2^4 + Z_3^4 = 0\}\) whose index is a generator of \(\KO_4(\ast)\), see~\cite[92]{LawsonMichelsohn:SpinGeometry}.

  We start with the case \(q = 0\) and prove the statement by induction on \(m\).
  Begin with \(m = 2\).
  Note that \(\KO_0^{\Z/2}(\ast) \otimes \Q\) is \(2\)-dimensonal.
  The equivariant index of \(\Z/2 \times V\), where the action on \(V\) is trivial, corresponds to the left-regular representation of \(\Z/2\).
  The other generator is obtained as follows.
  Consider the action of \(\Z/2\) by inversion \((z_1,z_2,z_3,z_4) \mapsto (\bar{z}_1, \bar{z}_2, \bar{z}_3, \bar{z}_4)\) on the standard \(4\)-torus \(\Torus^4 = \{(z_1,z_2,z_3,z_4) \in \C^4 \mid |z_i| = 1\}\).
  This action is free outside its \(16\) fixed points.
  Using the trivialization of the tangent bundle via the Lie group structure, one proves that this lifts to an action on the spinor bundle.
  The character corresponding to the equivariant index of this manifold can be computed using the Atiyah--Bott fixed point formula \cite[Theorem~8.35]{AB68:LefschetzApplications}.
  The signs that appear in the formula must be the same for each fixed point because the lift of the action to the spinor bundle was defined using a global trivialization.
  Hence the character satisfies \(\chi([1]) = -4 \csc\left(\frac{\pi}{2}\right) = -4 \neq 0\).
  Since it is non-trivial at the generator of \(\Z/2\), the representation generates \(\KO_0^{\Z/2}(\ast) \otimes \Q\) modulo the left-regular representation.

  Now let \(m > 2\).
  The induction hypothesis implies that the subspace \(\mathcal{I}_0 \subseteq \KO_0^{\Z/m}(\ast) \otimes \Q\) generated by those representations that are induced from a proper subgroup is generated by equivariant indices of \(4\)-dimensional closed spin \(\Z/m\)-manifolds, where the action is free outside a finite set of points.
  To obtain the remaining part of \(\KO_0^{\Z/m}(\ast) \otimes \Q\), we construct explicit examples.
  Set \(d \coloneqq \frac{m-1}{2}\) if \(m\) is odd and \(d \coloneqq \frac{m-2}{2}\) if \(m\) is even.
  For each \(k \in \{1, \dotsc, d\}\) that is coprime to \(m\) consider the action \(\alpha_k\) of \(\Z/m\) on the \(4\)-disk \(\disk^4 \subset \C^2\) such that the generator acts by scalar multiplication with \(\eu^{\frac{2 \pi \iu k}{m}} \in \sphere^1 \subset \C\), that is, $(z_1, z_2) \mapsto (\eu^{\frac{2 \pi \iu k}{m}} \cdot z_1, \eu^{\frac{2 \pi \iu k}{m}}\cdot  z_2) $.
  It is an orientation-preserving isometric action on the standard disk which, by coprimality, is free outside the fixed point at the origin.
  Moreover, the element corresponding to scalar multiplication on \(\C^2\) with \(\eu^{\frac{2 \pi \iu k}{m}}\) in \(\SO(4)\) has a lift of order \(m\) to the double covering \(\Spin(4)\).
  To see this, observe that the map \(\Spin(2) \to \SO(2)\) is the double self-covering \(\Sphere^1 \to \Sphere^1\).
  Thus there exists a lift \(\epsilon \in \Spin(2)\) of \(\eu^{\frac{2 \pi \iu k}{m}} \in \SO(2)\) such that \(\epsilon^m = (-1)^k\).
  Moreover, recall that \(\Spin(n)\) can be identified with a subgroup of the units in the even part of the Clifford algebra \(\Cl_n\) and that there are canonical identifications \(\Cl_{n_1} \tensgr \Cl_{n_2} = \Cl_{n_1 + n_2}\), see e.g.~\cite[Chapter~1]{LawsonMichelsohn:SpinGeometry}.
  Then \(\epsilon \tensgr \epsilon \in \Cl_2 \tensgr \Cl_2 = \Cl_4\) represents an element in \(\Spin(4)\) of order \(m\) which lifts the element in \(\SO(4)\) represented by multiplication with \(\eu^{\frac{2 \pi \iu k}{m}}\).
  Note that in general the lift to \(\Spin(2)\) does not have order \(m\) and this is why we need to work in dimension four here.
  In conclusion, this implies that \(\alpha_k\) extends to a spin action of \(\Z/m\) on \(\disk^4\).
  In particular, the lens space space obtained as the quotient manifold of the boundary sphere \(\sphere^3 / \alpha_k\) is a spin manifold with fundamental group \(\Z/m\).
  Since the spin bordism group \(\SpinBordism{3}{\BF \Z/m}\) is a torsion group, there is an \(n \in \N\) such that \(\bigsqcup_{i=1}^{n} \sphere^3 / \alpha_k\) is spin null-bordant over \(\BF \Z/m\) for each \(k\).
  Thus there exists a compact spin manifold \(W_k\) together with a free spin action \(\beta_k\) of \(\Z/m\) such that \(\partial(W_k, \beta_k) = \bigsqcup_{i=1}^{n} (\sphere^3, \alpha_k)\).
  Gluing \((W_k, \beta_k)\) along the boundary to \(\bigsqcup_{i=1}^{n} (\disk^3, \alpha_k)\), we obtain a spin manifold \(M_k\) together with a spin \(\Z/m\)-action which is free outside the origins of the copies of the disks.
  We will prove that the equivariant indices of \(M_k\) with \(\gcd(k, m) = 1\) together with the subspace \(\mathcal{I}_0\) generate \(\KO_0^{\Z/m}(\ast) \tens \Q\) and thereby complete the induction step.
  Recall that \(\KO_0^{\Z/m}(\ast) \tens \Q\) can be identified with the group of complex virtual representations whose characters are symmetric.
  Thus, by character theory, it suffices to establish the following claim:
  \begin{lem}\label{lem:Gershgorin}
    The vector space of symmetric functions \(\Fin^0(\Z/m)\) is generated by characters corresponding to the representations generating \(\mathcal{I}_0\) and those corresponding to the equivariant indices of \(M_k\), where \(k\) runs through all elements of \(\{1, \dotsc, d\}\) satisfying \(\gcd(k,m) =1\).
  \end{lem}
  \begin{proof}[Proof of \cref{lem:Gershgorin}]
  The dimension of \(\Fin^0(\Z/m)\) is equal to \(d+1\) if \(m\) is odd and \(d+2\) if \(m\) is even.
  It suffices to prove the statement modulo the characters which are supported on \([0], [m/2] \in \Z/m\) if \(m\) is even, and on \([0] \in \Z/m\) if \(m\) is odd, respectively. 
  This is because the representations corresponding to such characters lie in \(\mathcal{I}_0\) anyway.
  The remaining characters corresponding to \(\mathcal{I}_0\) are generated by
  \[
    \chi_k([l]) \coloneqq \begin{cases}
      1 & [l] = [\pm k],\\
      0 & \text{otherwise},
    \end{cases}  \quad k \in \{1, \dotsc, d\}, \gcd(k, m) \neq 1, \quad [l] \in \Z/m.
  \]
  To deal with the equivariant indices, we use the Atiyah--Bott fixed point formula \cite[Theorem~8.35]{AB68:LefschetzApplications} again.
  It implies that the character corresponding to \(M_k\), where \(\gcd(k,m) =1\), is given by
  \[
    \chi_k([l]) \coloneqq - \frac{n}{4} \csc^2 \left( \frac{\pi k l}{m} \right), \quad [l] \in \Z/m, [l] \neq 0.
  \]
  The factor \(n\) arises because by construction there are \(n\) identical fixed points in \(M_k\).
  By symmetry, the values on \([1],[2], \dotsc, [d] \in \Z/m\) determine such a character uniquely modulo the entries on \([0],[m/2]\in\Z/m\) (we do not consider the latter if \(m\) is odd).
  Thus, to prove the lemma, it suffices to establish that the matrix
  \[
    A = \begin{pmatrix}
    \chi_1([1]) & \chi_1([2]) & \hdots & \chi_1([d]) \\
    \chi_2([1]) & \chi_2([2]) &\hdots & \chi_2([d]) \\
    \vdots & & \ddots & \vdots \\
    \chi_d([1]) & \chi_d([2]) &\hdots & \chi_d([d])
  \end{pmatrix}
  \]
  is invertible.
  To that end, first consider \(k \in \{1, \dotsc, d\}\) with \(\gcd(k,m)=1\).
  Then the \(k\)-th row of the matrix \(A\) consists of a permutation of the entries of the first row.
  Recall that the first row reads explicitly as follows:
  \[
  -\frac{n}{4} \begin{pmatrix}\csc^2 \left( \frac{\pi}{m} \right) & \csc^2 \left( \frac{2\pi}{m} \right) & \hdots & \csc^2 \left( \frac{d \pi}{m} \right) \end{pmatrix}.
  \]
  The entry in the \(k\)-th row corresponding to \(\csc^2 \left( \frac{\pi}{m} \right)\) occurs in the \(l\)-th column, where \(l \in \{1, \dotsc, d\}\) is such that \(kl \equiv \pm 1 \mod m\).
  From this we can deduce that the entry \(-\frac{n}{4} \csc^2 \left( \frac{\pi}{m} \right)\) occurs precisely once in each of the columns of \(A\) whose column index \(l\) satisfies \(\gcd(l, m) = 1\).
  Thus, we can permute the rows of the matrix \(A\) in such a way that the \(k\)-th row has \(-\frac{n}{4} \csc^2 \left( \frac{\pi}{m} \right)\) on the diagonal position in the case of \(\gcd(k, m) = 1\) and is equal to the \(k\)-th unit vector in the case of \(\gcd(k, m) \neq 1\).
  Denote the matrix obtained from this row permutation by \(\tilde{A}\).

  The Gershgorin circle theorem implies that a matrix is invertible if for each row the sum of the absolute values of the non-diagonal entries is strictly smaller than the absolute value of the diagonal entry.
  We will verify this condition for the matrix \(\tilde{A}\).\footnote{The idea for this argument goes back to a comment on the MathOverflow question~\cite{OverflowCosecants}.}
  Therefore we will see that \(\tilde{A}\) and subsequently \(A\) are invertible, thereby finishing the proof of the lemma.
  By construction of \(\tilde{A}\), the rows of index \(k\) with \(\gcd(k, m) \neq 1\) are unit vectors with \(1\) on the diagonal position.
  So, the condition holds trivially in this case.
  In the case of \(\gcd(k, m) = 1\), the Gershgorin condition for the \(k\)-th row of \(\tilde{A}\) is equivalent to the inequality
  \begin{equation}
    \sum_{l=2}^{d} \csc^2 \left( \frac{l \pi}{m} \right) < \csc^2\left(\frac{\pi}{m} \right).
    \label{eq:GershgorinInequ}
  \end{equation}
  To prove \labelcref{eq:GershgorinInequ}, we make use of the elementary formula
   \begin{equation}\label{eq:CscFormula}
     \sum_{l=1}^{m-1} \csc^2\left(\frac{l \pi}{m}\right) = \frac{m^2-1}{3}.
   \end{equation}
   This formula first appeared as a step in an elementary solution to the \emph{Basel problem}, see~\cite[p.~56--57]{AZ01:TheBook} for this, where \labelcref{eq:CscFormula} is proved for \(m\) odd.
   Moreover, angle sum identities imply that \(\csc(x) + \csc(\pi/2+x) = 4 \csc(2x)\) for \(0 < x < \pi/2\).
   If we denote the left-hand side of \labelcref{eq:CscFormula} by \(S_m\), then this implies the recurrence relation \(S_{2m} = 4 S_{m} + 1\).
   As a consequence of this, \labelcref{eq:CscFormula} also holds if \(m\) is even.
   
   This implies
  \[
    \sum_{l=1}^{d} \csc^2\left(\frac{l \pi}{m}\right) = \begin{cases}
    \frac{m^2-1}{6} & \text{\(m\) odd} \\
    \frac{m^2-4}{6} & \text{\(m\) even}
  \end{cases} < \frac{m^2}{6}.
  \]
  Thus
  \begin{align*}
    \sum_{l=2}^{d} \csc^2 \left( \frac{l \pi}{m} \right) &< \frac{m^2}{6} - \csc^2\left(\frac{\pi}{m} \right)  \\
    &< 2\frac{m^2}{\pi^2} - \csc^2\left(\frac{\pi}{m} \right) \\
    &\leq 2\csc^2\left(\frac{\pi}{m} \right) - \csc^2\left(\frac{\pi}{m} \right) = \csc^2\left(\frac{\pi}{m} \right),
  \end{align*}
  where we used the elementary estimates \(\pi^2/2 < 6\) and \(x^{-2} \leq \csc^2(x)\) for \(x \in \R \setminus \{0\}\).
  This establishes \labelcref{eq:GershgorinInequ} as required.
\end{proof}
The lemma completes the induction step and hence the proof of \cref{prop:finCyc} in the case \(q = 0\).

Continue with the proof in the case \(q =1\).
We again set \(d \coloneqq \frac{m-1}{2}\) if \(m\) is odd and \(d \coloneqq \frac{m-2}{2}\) if \(m\) is even.
We will use the following lemma.
\begin{lem}\label{lem:2DExamples}
  Suppose that \(m > 2\) and choose \(k \in \{1, \dotsc, d\}\).
  There exists a 2-dimensional spin \(\Z/m\)-manifold such that the character \(\chi\) corresponding to its equivariant index satisfies \(\chi([\pm k]) \neq 0\).
\end{lem}
\begin{proof}[Proof of \cref{lem:2DExamples}]
  The order of the subgroup generated by \([k]\) is also strictly larger than \(2\) because \(d < \frac{m}{2}\).
  Thus we can assume without loss of generality that \([k]\) generates \(\Z/m\) because otherwise we can induce up from the subgroup generated by it.
  
  As a further preparation, we claim that if there exists a surjective homomorphism \(\varphi \colon \Z/{m'} \twoheadrightarrow \Z/{m}\) and the lemma is proved for \(m\) and all generators \([k]\) of \(\Z/m\), then the lemma also holds for \(m'\) and all generators \([k']\) of \(\Z/m'\).
  This is because any such \(\Z/{m}\)-manifold with corresponding index character \(\chi\) can be turned into a \(\Z/m'\)-manifold with index character \(\chi' = \chi \circ \varphi\) via the homomorphism \(\varphi\).
  Since \(\varphi\) is surjective, it maps each generator of \(\Z/m'\) to a generator of \(\Z/m\) and the claim follows.
  
  Now, since \(m > 2\), there exists a surjective homomorphism \(\Z/m \twoheadrightarrow \Z/4\) or a surjective homomorphism \(\Z/m \twoheadrightarrow \Z/r\), where \(r\) is an odd number.
  Thus by these preparations it suffices to prove the lemma for the cases \(m = 4\) and \(m\) odd.

  Start with \(m = 4\).
  The action of \(\Z/2\) by inversion on the standard \(2\)-torus lifts to an action of \(\Z/4\) on the spinor bundle which is non-trivial on the generator.
  The Atiyah--Bott fixed point formula~\cite[Theorem~8.35]{AB68:LefschetzApplications} shows that the corresponding character satisfies \(\chi([k]) = 2 \iu \csc\left(\frac{\pi}{2}\right) = 2 \iu \neq 0\) if \([k]\) is a generator of \(\Z/4\).

  To deal with an arbitrary odd \(m\), we use an analogous construction as in the induction step in the case \(q=0\).
  That is, consider the action \(\alpha\) of \(\Z/m\) on the \(2\)-disk \(\disk^2 \subset \C\) such that the generator acts by multiplication with \(\eu^{\frac{2 \pi \iu}{m}}\).
  It is an orientation-preserving isometric action which is free outside the fixed point at the origin.
  Moreover, \(\eu^{\frac{2 \pi \iu}{m}} \in \SO(2)\) has a lift of order \(m\) to \(\Spin(2)\).
  Note that here we used that \(m\) is odd because in dimension two this fails for even \(m\).
  Thus this construction yields a spin action on the \(2\)-disk.
  Furthermore, by the same argument as before, some number \(n \in \N\) of copies of the boundary circle is equivariantly null-bordant with a free action.
  Gluing these null-bordisms with the disk, we obtain a \(2\)-dimensional spin \(\Z/m\)-manifold \(M\) which is free outside \(n\) identical fixed point.
  Now the Atiyah--Bott fixed point formula~\cite[Theorem~8.35]{AB68:LefschetzApplications} implies that the character \(\chi\) corresponding to the equivariant index of \(M\) for a generator \([k] \in \Z/m\) satisfies
  \[
    \chi([k]) = \frac{n \iu}{2} \csc\left(\frac{k \pi}{m} \right) \neq 0 \qedhere
  \]
\end{proof}
To see how \cref{lem:2DExamples} implies the claim of \cref{prop:finCyc} in the case \(q =1\), recall that \(\KO_2^{\Z/m}(\ast) \otimes \Q\) can be identified with complex virtual representations whose characters are antisymmetric.
Since \(\KO_2^{\Z/2}(\ast) \otimes \Q = 0\), we can assume that \(m > 2\).
The space of antisymmetric functions \(\Fin^1(\Z/m)\) has dimension \(d\) and is generated by
\[
  \chi_k([l]) = \begin{cases}
  \pm1 & [l] = [\pm k], \\
  0 & \text{otherwise},
\end{cases} \quad
k \in \{1, \dotsc, d\}, \quad [l] \in \Z/m.
\]
Fix \(k \in \{1, \dotsc, d\}\).
It follows from the already established case \(q =0\) that there exist \(4\)-dimensional spin \(\Z/m\)-manifolds which are free outside a finite set of fixed points whose corresponding characters admit a linear combination which is equal to
\[
  \psi_k([l]) = \begin{cases}
  1 & [l] = [\pm k], \\
  0 & \text{otherwise},
\end{cases} \quad [l] \in \Z/m.
\]
Taking products of these with an example given by \cref{lem:2DExamples}, yields \(6\)-dimensional spin \(\Z/m\)-manifolds whose corresponding characters admit a linear combination which is equal to a multiple of \(\chi_k\).
Moreover, these products are free outside a submanifold of codimension \(4\) with trivial normal bundle because they are products of \(4\)-dimensional \(\Z/m\)-manifolds, which are free outside a submanifold of codimension \(4\), with arbitrary \(2\)-dimensional \(\Z/m\)-manifolds.
Since the functions \(\chi_k\) with \(k\) running through \(\{1, \dotsc, d\}\) generate \(\Fin^1(\Z/m)\), we deduce by character theory that \(\KO_2^{\Z/m}(\ast) \otimes \Q\) is generated by examples of these type.
This concludes the proof of the case \(q = 1\) and thus the proof of \cref{prop:finCyc}.
\end{proof}

In the proof of our main result, we will need that the manifolds from \cref{prop:finCyc} can be endowed with an invariant metric of positive scalar curvature in a neighborhood around the submanifold where the action is non-free.
This can be achieved by way of \emph{torpedo metrics}:

\begin{defi}\label{defi:torpedo}
  Let $\disk^d$ be the $d$-disk.
  A \emph{torpedo metric} is an $\OG(d)$-invariant positive scalar curvature metric $g_{\mathrm{tor}} \in \Riem^+(\disk^d)$ which is of product structure near the boundary $\bd \disk^d = \sphere^{d-1}$, agrees with the round metric at the boundary, and near the origin $0 \in \disk^d$ agrees with the round hemispheric metric on the disk.
\end{defi}
Torpedo metrics are a standard tool in the study of the topology of positive scalar curvature.
They always exists if $d \geq 3$; for a detailed construction see for instance Walsh~\cite[Chapter 1]{Walsh11pscGenMorse}.
We use them in the following:
\begin{prop}\label{prop:tubularPSC}
  Let \(H\) be a finite group acting on a closed manifold \(M\) and \(N \subset M\) an \(H\)-invariant submanifold of codimension \(d \geq 3\).
  Then an \(H\)-invariant tubular neighborhood of \(N\) can be endowed with an \(H\)-invariant metric of positive scalar curvature which is collared near the boundary.
\end{prop}
\begin{proof}
  Let \(\nu \to N\) denote the normal bundle of \(N\) in \(M\) and start with any \(H\)-equivariant tubular neighborhood embedding \(\nu \hookrightarrow M\).
  Fix an \(H\)-invariant Riemannian metric on \(N\) and an \(H\)-invariant fiberwise Euclidean metric on the normal bundle \(\nu\) and form the associated unit disk bundle \(\disk \nu \to N\).
  Using the torpedo metric fiberwise (this makes sense because of \(\OG(d)\)-invariance), we obtain an \(H\)-invariant fiberwise Riemannian metric on \(\disk \nu \to N\) which is fiberwise collared near the boundary.
  After choosing an \(H\)-invariant horizontal bundle for \(\disk \nu \to N\), we can use these fiberwise metrics to construct an \(H\)-invariant Riemannian metric on \(\disk \nu\), collared near the boundary, which turns the map \(\disk \nu \to N\) into a Riemannian submersion.
  Applying O'Neill's formulas~\cite[Chapter~9~D]{Besse08:Einstein} after—if necessary—shrinking the fibers, we see that this metric has positive scalar curvature.
\end{proof}

\subsection{Main results}\label{subsec:mainResults}
Finally, we prove our main theorem and state its corollaries.
\begin{thm}
  Let $p \in \{0,1,2,3\}$ and $k \geq 1$.
  Then the image of the composition
  \begin{equation*}
    \StolzPartialPSC[\Gamma]{4k+p}{\EFin \Gamma, \EF \Gamma} \tens \C \xrightarrow{\omega \tens \C} \SpinBordism[\Gamma]{4k + p}{\EFin \Gamma} \tens \C \xrightarrow{\beta \tens \C} \KO_{p}^\Gamma(\EFin \Gamma) \tens \C
  \end{equation*}
  contains
  \begin{equation*}
    \bigoplus_{l < k} \HZ_{4l + p}(\Gamma; \Fin^0\Gamma) \oplus \HZ_{4l - 2 + p}(\Gamma; \Fin^1 \Gamma)
  \end{equation*}
  with respect to the equivariant delocalized Pontryagin character.
\end{thm}
\begin{proof}
  By \citeauthor{matthey:delocChern}~\cite[Theorem~1.3 and Section~7]{matthey:delocChern}, the \emph{complex} K-homology group $\KU_{p}^\Gamma(\EFin \Gamma) \tens \C$ is generated by the images of maps of the form
  \begin{align*}
    \KU_{p}^\Lambda(\EF \Lambda) \tens \KU^{H}_{0}(\pt) \tens \C &\to \KU_p^\Gamma(\EFin \Gamma) \tens \C, \\
    x \tens y  &\mapsto \psi_\ast(x \times y),
  \end{align*}
  where $\Lambda \subseteq \Gamma$ is a subgroup, $\gamma \in \Gamma$ is an element of some finite order $m$ which commutes with $\Lambda$, $H = \langle \gamma \rangle \iso \Z/m$, and $\psi \colon \Lambda \times H \to  \Gamma$ is the homomorphism $(\lambda, \gamma^l) \mapsto \lambda \gamma^l$.
  
  As a consequence, using the isomorphism \(\KU_{\ast}^{\bullet} \tens \C \cong (\KO_{\ast}^{\bullet} \oplus \KO_{\ast+2}^{\bullet})\tens \C\) of equivariant homology theories (see~\cite[Proposition~2.1]{BZ17LowDegree}) on both sides and sorting the summands appropriately, we obtain that $\KO_{p}^\Gamma(\EFin \Gamma) \tens \C$ is generated by the images of maps of the form
  \begin{align*}
    \KO_{p - 2q}^\Lambda(\EF \Lambda) \tens \KO^{H}_{2q}(\pt) \tens \C &\to \KO_p^\Gamma(\EFin \Gamma) \tens \C, \\
    x \tens y  &\mapsto \psi_\ast(x \times y),
  \end{align*}
  where $q \in \{0,1\}$ and, as above, $\Lambda \subseteq \Gamma$ is a subgroup, $\gamma \in \Gamma$ is an element of some finite order $m$ which commutes with $\Lambda$, $H = \langle \gamma \rangle \iso \Z/m$, $\psi \colon \Lambda \times H \to  \Gamma$ is the homomorphism $(\lambda, \gamma^l) \mapsto \lambda \gamma^l$.
  Moreover, it follows from the construction of Matthey's handicrafted Chern character~\cite[Theorem~1.4]{matthey:delocChern} that there is a commutative diagram
  \[
  \begin{tikzcd}
  \KO^\Lambda_{p-2q}(\EF \Lambda) \tens \KO^H_{2q}(\ast) \tens \C \rar{\psi_\ast \circ \times} \dar{\ph \tens \phub_H}
    & \KO_p^\Gamma(\EFin \Gamma) \tens \C \dar{\phub_\Gamma} \\
  \bigoplus_{k\in \Z} \HZ_{p-2q+4k}(\Lambda;\Q) \tens \HZ_0(H, \Fin^q H) \rar{\psi_\ast \circ \times}
    & \bigoplus_{k \in \Z} \HZ_{p+4k}(\Gamma; \Fin^q \Gamma).
  \end{tikzcd}
    \]
    To be precise, we again first deduce the complex variant of this from~\cite[Theorem~1.4]{matthey:delocChern} which, in turn, restrict to the real version via the isomorphism \(\KU_{\ast}^{\bullet} \tens \C \cong (\KO_{\ast}^{\bullet} \oplus \KO_{\ast+2}^{\bullet})\tens \C\); see also~\cite[Proof~of~Proposition~2.2]{BZ17LowDegree}.
    From this diagram we then deduce that if $x \in \KO_{p-2q}^\Lambda(\EF \Lambda) \otimes \Q$ with $\ph(x) \in \HZ_{p - 2q + 4l}(\Lambda; \Q)$ for some $l \in \Z$ and $y \in \KO_{2q}^{\Z/m}(\pt)$, then
  \begin{equation*}
    \phub_\Gamma \left(\psi_\ast(x \times y) \right) \in \HZ_{p+4l-2q}(\Gamma; \Fin^q \Gamma).
  \end{equation*}
  Since the (non-equivariant) Pontryagin character for \(\BF \Lambda\) and the equivariant Pontryagin character for \(H = \Z/m\) are both isomorphisms, this in particular means that the group $\HZ_{p+4l-2q}(\Gamma; \Fin^q \Gamma)$ is generated by such elements.

  Hence it suffices to show that for all $\Lambda \subseteq \Gamma$, $\gamma \in \Gamma$ as in the previous paragraph and each $q \in \{0,1\}$, $l < k$, the image of $\beta \circ \omega$ contains all $\psi_\ast(x \times y)$ with $\ph(x) \in \HZ_{p-2q+4l}(\Lambda; \Q)$ and $y \in \KO_{2q}^H(\pt)$.

  By \cref{prop:framedBordism}, we may assume that there is a closed spin manifold $M$ of dimension $(p-2q+4l)$ and a map $\sigma \colon M \to \B \Lambda$  with $\bar{\sigma}_\ast [\bar{M}]_\KO^\Lambda = x$, where $\bar{M} \to M$ is the $\Lambda$-covering induced by $\sigma$.
  Furthermore, by \cref{prop:finCyc}, we can assume that there exists a closed $H = \Z/m $-manifold $N$ of dimension $4(k-l) + 2q$, which is free outside a submanifold \(N_{0} \subset N\) of codimension \(4\), such that $y = p_\ast [N]^H_\KO$, where $p \colon N \to \pt$.

Let \(N_+ \subset N\) be an \(H\)-invariant tubular neighborhood of \(N_0\) that is endowed with an \(H\)-invariant Riemannian metric of positive scalar curvature \(g \in \Rm^+(N_+)\) collared near the boundary (see \cref{prop:tubularPSC}).
Define $N_- \subset N$ to be the complement of the interior of $N_+$.
Then $N_-$ is also $H$-invariant and $N$ is the union of $N_-$ and $N_+$.
 The $H$-action on $N_{-}$ is free because all non-free points are contained in $N_0$.
  Thus there is an $H$-equivariant classifying map $\nu_{-} \colon N_{-} \to \EF H$.
  Since $N_+$ is a disk bundle on which $H$ acts by orthogonal transformations between the fibers, we can extend the map $\nu_{-}|_{\bd N_-}$ to an $H$-equivariant map $\nu_{+} \colon N_{+} \to \operatorname{Cone}(\EF H)$ using the fiberwise radius as the cone coordinate.
  Here $\operatorname{Cone}(\EF H)$ denotes the cone over $\EF H$.
  Since \(H\) is finite, it is a model for the classifying space $\EFin H$, and we will use this instead of the one-point space in the following.
Next we choose a $\Lambda$-invariant Riemannian metric $g_{\bar{M}}$ on $\bar{M}$ such that $g_{\bar{M}} \oplus g$ has uniformly positive scalar curvature (by first finding an arbitrary $\Lambda$-invariant metric on $\bar{M}$ and rescaling it if necessary).

  Let $W_{\pm}$ be the cocompact \((4k+p)\)-dimensional $\Lambda \times H$-spin manifold defined as $\bar{M} \times N_{\pm}$ and let
  \begin{equation*}
    \xi := [(W_{-}, W_{+}, \bar{\sigma} \times \nu_{\pm}, g_{\bar{M}}  \oplus g)] \in \StolzPartialPSC[\Lambda \times H]{4k+p}{\EF \Lambda \times \EFin H, \EF \Lambda \times \EF H}
  \end{equation*}
  If we apply the forgetful map $\omega$ from \cref{rem:partialPSCForget}, we obtain
  \begin{equation*}
    \omega(\xi) = \bar{\sigma}_\ast[\bar{M}]^\Lambda \times \nu_\ast[N]^H\in \SpinBordism[\Lambda \times H]{4k+p}{\EF \Lambda \times \EFin H},
  \end{equation*}
  where $\nu := \nu_- \cup \nu_+ \colon N = N_- \cup N_+ \to \EFin H = \Cone(\EF H)$.
  Thus
  \begin{equation}
    (\beta \circ \omega)(\xi) = \bar{\sigma}_\ast[\bar{M}]_\KO^\Lambda \times \nu_\ast[N]_\KO^H = x \times y \in \KO^{\Lambda \times H}_p(\EF \Lambda \times \EFin H) \tens \C.\label{eq:liftProductxy}
    \end{equation}
  Now consider the induction homomorphism (compare \cref{rem:induction})
  \begin{equation*}
  \resizebox{\linewidth}{!}{$
    \induction_\psi \colon \StolzPartialPSC[\Lambda \times H]{4k+p}{\EF \Lambda \times \EFin H, \EF \Lambda \times \EF H}
    \to  \StolzPartialPSC[\Gamma]{4k+p}{\Gamma \times_\psi (\EF \Lambda \times \EFin H, \EF \Lambda \times \EF H)}.$}
  \end{equation*}
  Since the $\Gamma$-action is proper on $\Gamma \times_\psi (\EF \Lambda \times \EFin H)$ and free on $\Gamma \times_\psi (\EF \Lambda \times \EF H)$,
  there is a map $f \colon \Gamma \times_\psi (\EF \Lambda \times \EFin H, \EF \Lambda \times \EF H) \to (\EFin \Gamma, \EF \Gamma)$.
  Then we set
  \begin{equation*}
    \Xi := f_\ast \induction_\psi (\xi)  \in \StolzPartialPSC[\Gamma]{4k+p}{\EFin \Gamma, \EF \Gamma}
  \end{equation*}
  Finally, \labelcref{eq:liftProductxy} implies that $(\beta \circ \omega)(\Xi) = \psi_\ast(x \times y) \in \KO_p^\Gamma(\EFin \Gamma)$.
\end{proof}

\begin{cor}
  Let $n \geq 4$.
  If the rational homological dimension of $\Gamma$ is at most $n -3$, then the composition
  \begin{equation*}
    \StolzPartialPSC[\Gamma]{n}{\EFin \Gamma, \EF \Gamma} \xrightarrow{\omega} \SpinBordism[\Gamma]{n}{\EFin \Gamma} \to \KO_n^\Gamma(\EFin \Gamma)
  \end{equation*}
  is rationally surjective.
\end{cor}

\begin{cor}\label{cor:imageRelIndex}
  Let $p \in \{0,1,2,3\}$ and $k \geq 1$.
  Then the image of the relative index map $\alpha \tens \C \colon \StolzRel{4k + p}{\BF \Gamma} \tens \C \to \KO_{4k + p}(\CstarRed \Gamma) \tens \C$
  contains
  \begin{equation*}
    \mu \left( (\phub_\Gamma)^{-1}\left(\HZ_{4l -2q + p}(\Gamma;  \Fin^q\Gamma)\right)\right) \subseteq \KO_{p}(\CstarRed \Gamma) \tens \C
  \end{equation*}
  for each $l < k$ and $q \in \{0,1\}$.
\end{cor}

\begin{cor}
      Let $p \in \{0,1,2,3\}$ and $k \geq 1$.
    Suppose that the Baum--Connes assembly map $\mu \colon \KO_\ast^\Gamma(\EFin \Gamma) \to \KO_\ast(\CstarRed \Gamma)$ is rationally injective.
    Then the rank of $\StolzRel{4k + p}{\BF \Gamma}$ is at least the dimension of
    \begin{equation*}
      \bigoplus_{l < k} \bigoplus_{q \in \{0,1\}} \HZ_{4l -2q + p}(\Gamma;  \Fin^q\Gamma).
    \end{equation*}
\end{cor}

\begin{cor}\label{cor:relIndexSurj}
  Let $n \geq 4$ and the rational homological dimension of $\Gamma$ be at most $n -3$. Suppose that the real Baum--Connes assembly map for $\Gamma$ is rationally surjective.
  Then the relative index map $\alpha \colon \StolzRel{n}{\BF \Gamma} \to \KO_{n}(\CstarRed \Gamma)$ is rationally surjective.
\end{cor}

We also obtain consequences for the higher rho-invariant by simply applying the boundary maps in Stolz' positive scalar curvature sequence and the analytic surgery sequence of Higson and Roe, respectively.
Thus we need to identify the part in equivariant group homology that goes to zero under the boundary map.
To that end, let $\Fin^q_0 \Gamma := \{ f \in \Fin^q_0 \Gamma \mid f(1) = 0\}$.
Then $\Fin^0 \Gamma = \C \oplus \Fin^0_0 \Gamma$ and $\Fin^1 \Gamma = \Fin^1_0 \Gamma$.
\begin{lem}
  If the Baum--Connes assembly map is rationally injective, then $\bd \circ \mu \circ \phub_\Gamma^{-1}$ maps
  \begin{equation*}
    \bigoplus_{l \in \Z} \bigoplus_{q \in \{0,1\}} \HZ_{4l - 2q + p}(\Gamma; \Fin^q_0 \Gamma)
  \end{equation*}
  injectively into the structure group $\SG_{p-1}^\Gamma(\EF \Gamma) \tens \C$.
\end{lem}
\begin{proof}
By construction of $\Fin^q_0$ ($q \in \{0,1\}$), we have split short exact sequences $0 \to \HZ_\ast(\Gamma; \C) \to \HZ_\ast(\Gamma; \Fin^0 \Gamma) \to \HZ_\ast(\Gamma; \Fin^0_0 \Gamma) \to 0$ and $0 \to 0 \to \HZ_\ast(\Gamma; \Fin^1 \Gamma) \to \HZ_\ast(\Gamma; \Fin^1_0 \Gamma) \to 0$.
From this we obtain the following commutative diagram, where the top and bottom rows are exact:
  \begin{equation*}
  \resizebox{\linewidth}{!}{
    \begin{tikzcd}[column sep=small, ampersand replacement=\&]
      \bigoplus_{l \in \Z} \HZ_{4l + p}(\Gamma; \C) \rar[hook]
        \& \bigoplus_{l \in \Z} \bigoplus_{q \in \{0,1\}} \HZ_{4l - 2q + p}(\Gamma; \Fin^q \Gamma) \rar[two heads] \lar[two heads, bend right=20]
        \& \bigoplus_{l \in \Z} \bigoplus_{q \in \{0,1\}} \HZ_{4l - 2q + p}(\Gamma; \Fin^q_0 \Gamma) \lar[bend right=20, hook] \ar[dd, hook, dashed, "\bd \circ \mu \circ \phub_\Gamma^{-1}"]
        \\
         \& \KO_p^\Gamma(\EFin \Gamma) \tens \C \ar[u, "\phub_\Gamma" swap, "\cong"] \dar[hook, "\mu \tens \C"]\\
         \KO_p(\BF \Gamma) \tens \C \ar[ur] \ar[uu, "\ph_\Gamma" swap, "\cong"] \rar["\nu \tens \C"] \& \KO_p(\CstarRed \Gamma) \tens \C \rar["\bd \tens \C"]
          \& \SG_{p-1}^\Gamma(\EF \Gamma) \tens \C
    \end{tikzcd}}
  \end{equation*}
  A straightforward diagram chase implies that the vertical arrow on the right-hand side must be injective.
\end{proof}

\begin{cor}
      Let $p \in \{0,1,2,3\}$ and $k \geq 1$.
    Suppose that the Baum--Connes assembly map $\mu \colon \KO_\ast^\Gamma(\EFin \Gamma) \to \KO_\ast(\CstarRed \Gamma)$ is rationally injective.
    Then the rank of $\StolzPos{4k + p - 1}{\BF \Gamma}$ is at least the dimension of
    \begin{equation*}
      \bigoplus_{l < k} \bigoplus_{q \in \{0,1\}} \HZ_{4l -2q + p}(\Gamma;  \Fin^q_0\Gamma).
    \end{equation*}
\end{cor}

\begin{cor}\label{cor:RhoSurj}
  Let $n \geq 4$ and the rational homological dimension of $\Gamma$ be at most $n -3$. Suppose that the real Baum--Connes assembly map for $\Gamma$ is a rational isomorphism.
  Then the higher rho-invariant \(\rho \colon \StolzPos{n-1}{\BF \Gamma} \to \SG_{n-1}^\Gamma(\EF \Gamma)\) is rationally surjective.
\end{cor}
\begin{proof}
  Rational injectivity of the Baum--Connes assembly map implies rational surjectivity of $\bd \colon \KO_n(\CstarRed \Gamma) \to \SG_{n-1}^\Gamma(\EF \Gamma)$.
  Thus the result follows from \cref{cor:relIndexSurj}.
\end{proof}

Finally, we draw conclusions for the Bott stabilized index maps (see \labelcref{eq:bottStabilize,eq:stabilizedRelIndex} in the introduction):

\begin{cor}\label{cor:relIndexStableSurj}
 Suppose that the real Baum--Connes assembly map for $\Gamma$ is rationally surjective.
  Then the stabilized relative index map
  \[\alpha[\Btmfd^{-1}] \colon \StolzRel{n}{\BF \Gamma}[\Btmfd^{-1}] \to \KO_{n}(\CstarRed \Gamma)\]
   is rationally surjective.
\end{cor}

\begin{proof}
  \Cref{cor:imageRelIndex} implies that the image of every equivariant homology class under the Baum--Connes assembly eventually lies in the image of $\alpha \otimes \C \colon \StolzRel{n}{\BF \Gamma} \to \KO_n(\CstarRed \Gamma)$ for sufficiently large $n$.
  Hence the image of every such class lies in the image of the stabilized relative index map.
  Thus surjectivity of the Baum--Connes assembly map implies that every class in $\KO_\ast(\CstarRed \Gamma)$ lies in the image of $\alpha[\Btmfd^{-1}]$.
\end{proof}

For the Bott stabilized higher rho-invariant we conclude as in \cref{cor:RhoSurj}:
\begin{cor}
 Suppose that the real Baum--Connes assembly map for $\Gamma$ is rationally bijective.
  Then the stabilized rho-invariant
  \[\rho[\Btmfd^{-1}] \colon \StolzPos{n-1}{\BF \Gamma}[\Btmfd^{-1}] \to \SG_{n-1}^\Gamma(\EF \Gamma)\]
   is rationally surjective.
\end{cor}

\section{The two definitions of the group of concordance classes}
\label{sec:twoDefsOfConcordanceGp}
Let $\tilde{\pi}_0(\Riem^+(M))$ be the set of concordance classes of positive scalar curvature metrics on a high-dimensional closed spin manifold $M$.
The group $\StolzRel{n+1}{\BF \Gamma}$ acts freely and transitively on $\tilde{\pi}_0(\Riem^+(M))$ \cite[Theorem~1.1]{Stolz98Concordance}.
After fixing a base-point $g_0 \in \Riem^+(M)$ this endows $\tilde{\pi}_0(\Riem^+(M))$ with a group structure.
Recently, \citeauthor{WY15FinitePart}~\cite{WY15FinitePart} also defined a group structure on $\tilde{\pi}_0(\Riem^+(M))$ depending on a base-point $g_0 \in \Riem^+(M)$ and denoted the resulting group by $P(M)$.
We show that these two group structures are the same.
\begin{prop}
  Let $M$ be a closed spin manifold of dimension $n \geq 5$ with $\pi_1 M = \Gamma$.
  Then the group structure induced by the free and transitive action of $\StolzRel{n+1}{\BF \Gamma}$ on $\tilde{\pi}_0(\Riem^+(M))$ is the same as in the group $P(M)$ of Weinberger and Yu provided that the same base-point $g_0 \in \Riem^+(M)$ is used in the definition of both group structures.
\end{prop}
\begin{proof}
  Start with the definition of the action of $\StolzRel{n+1}{\BF \Gamma}$ on $\tilde{\pi}_0(\Riem^+(M))$ following \citeauthor{Stolz98Concordance}~\cite[24]{Stolz98Concordance}:
  Let
  \begin{equation*}
    i \colon \tilde{\pi}_0(\Riem^+(M)) \times  \tilde{\pi}_0(\Riem^+(M)) \to \StolzRel{n+1}{\BF \Gamma}
  \end{equation*}
  be the map which takes a pair of metrics of positive scalar curvature $(g, g^\prime)$ to the element $i(g,g^\prime)$ of $\StolzRel{n+1}{\BF \Gamma}$ represented by $M \times [0,1]$ endowed with the metric $g \cup g^\prime$ on $\bd(M \times [0,1]) = M \times \{0\} \cup M \times \{1\}$.
  Then the action of $\StolzRel{n+1}{\BF \Gamma}$ on $\tilde{\pi}_0(\Riem^+(M))$ is implicitly defined by the requirement $i(g, g^\prime) \cdot [g] = [g^\prime]$ for all $g, g^\prime \in \Riem^+(M)$.
  Moreover,  $i(g, g^\prime) + i(g^\prime, g^{\prime\prime}) = i(g, g^{\prime\prime})$ for all $g, g^\prime, g^{\prime\prime} \in \Riem^+(M)$.

  The group $P(M)$ of Weinberger and Yu is defined as follows~\cite[2789]{WY15FinitePart}.
  Let $\Intv = [0,1]$.
  Consider the connected sum $(M \times \Intv) \sharp (M \times \Intv)$ which we perform away from the boundaries.
  Then apply further surgeries to $(M \times \Intv) \sharp (M \times \Intv)$ to remove the kernel of the homomorphism
  \[ \pi_1((M \times \Intv) \sharp (M \times \Intv)) = \Gamma \ast \Gamma \to \Gamma.\]
  The resulting manifold is the generalized connected sum $(M \times \Intv) \natural (M \times \Intv)$.
  By construction its fundamental group is $\Gamma$.
  Now let $g, g^\prime \in \Riem^+(M)$.
  The boundary of $(M \times \Intv) \natural (M \times \Intv)$ has four components, two being $M$ and the other two being $-M$.
  We put $g_0$ on one of the components $M$ and each of $g, g^\prime$ on one of the two components $-M$.
  Then the surgery theorem implies that there is a metric of positive scalar curvature $h$ on $(M \times \Intv) \natural (M \times \Intv)$ extending $g_0 \cup g \cup g^\prime$ with product structure near the boundary.
  Restricting $h$ to the remaining boundary component $M$ yields a metric of positive scalar curvature $g^{\prime \prime} \in \Riem^+(M)$.
  Finally, the group structure of Weinberger and Yu is defined by $[g] + [g^\prime] := [g^{\prime\prime}]$.

 To prove that the group structures agree, we need to show
  \[ i(g_0, g^{\prime\prime}) = i(g_0, g) + i(g_0, g^\prime).\]
  Indeed, let $x \in \StolzRel{n+1}{\BF \Gamma}$ be the element represented by $(M \times \Intv) \natural (M \times \Intv)$ endowed with the metric $g_0 \cup g \cup g^{\prime\prime} \cup g^\prime$ on its boundary $\bd ((M \times \Intv) \natural (M \times \Intv)) = M \cup -M \cup M \cup -M$.
  On the one hand, since the metric extends to the metric $h$ of positive scalar curvature on all of $(M \times \Intv) \natural (M \times \Intv)$, the element $x$ vanishes in $\StolzRel{n+1}{\BF \Gamma}$.
  On the other hand, undoing the surgeries and the connected sum shows that $(M \times \Intv) \natural (M \times \Intv)$ is bordant relative boundary to $M \times \Intv \sqcup M \times \Intv$.
  Hence $0 = x = i(g_0, g) + i(g^{\prime\prime}, g^\prime)$ and so $i(g_0, g) = i(g^\prime, g^{\prime\prime})$.
  We conclude that $i(g_0, g) + i(g_0, g^\prime) = i(g^\prime, g^{\prime\prime}) + i(g_0, g^\prime) = i(g_0, g^{\prime\prime})$.
\end{proof}

\section*{Acknowledgements}
Rudolf Zeidler gratefully acknowledges the hospitality of the Department of Mathematics at the Texas A\&M University, where this project was initiated. He also wishes to thank Johannes Ebert for several helpful conversations.
The authors would also like to thank the referee for their careful reading of the manuscript and for useful suggestions.

Guoliang Yu was partially supported by NSF 1564398, 1700021, 2000082, and the Simons Fellows Program.
Zhizhang Xie was partially supported by NSF 1500823, 1800737, and 1952693.
Rudolf Zeidler was partially supported by the Deutsche Forschungsgemeinschaft (DFG, German Research Foundation) through the CRC 878 and under Germany's Excellence Strategy EXC 2044-390685587, Mathematics Münster: Dynamics -- Geometry -- Structure in Münster.

 \printbibliography
 \vspace{2ex}
\end{document}